\documentclass[12pt]{amsart}
\usepackage{exscale,aliascnt,amssymb,amsthm,amsxtra,latexsym,mathrsfs,mathtools,paralist,varioref,verbatim}
\usepackage[T1]{fontenc}
\usepackage[margin=1in]{geometry}
\usepackage{tikz}

\usepackage[pagebackref]{hyperref} 

\usepackage[capitalize]{cleveref} 

\usepackage[abbrev,backrefs,msc-links]{amsrefs} 

\DefineSimpleKey{bib}{article-id}
\DefineSimpleKey{bib}{category}

\newcommand{\arxiv}[1]{\href{http://arxiv.org/abs/math/#1}{arXiv:#1 [\BibField{category}]}}

\BibSpec{arxiv}{%
    +{}  {\PrintAuthors}                {author}
    +{,} { \textit}                     {title}
    +{}  { \PrintDatePV}                {date}
    +{,} { \eprintpages}                {pages}
    +{,} { \arxiv}                      {article-id}
    +{.} { }                            {note}
    +{.} {}                             {transition}
}

\theoremstyle{plain}
\newtheorem{thm}{Theorem}

\newaliascnt{cor}{thm}

\aliascntresetthe{cor}
\newaliascnt{lem}{thm}

\aliascntresetthe{lem}
\newaliascnt{prop}{thm}

\aliascntresetthe{prop}

\theoremstyle{definition}
\newaliascnt{defn}{thm}
\newtheorem{defn}[defn]{Definition}
\aliascntresetthe{defn}

\theoremstyle{remark}

\newtheorem{ex}[thm]{Example}

\DeclareMathOperator{\convhull}{Conv}

\DeclareMathOperator{\newt}{Newt}

\DeclareMathOperator{\pull}{pull}

\DeclareMathOperator{\spec}{Spec}

\newcommand*{\monoid}{\mathsf{S}}

\newcommand*{\sheaf}[1]{\mathcal#1}

\newcommand*{\QQ}{\mathbb{Q}}
\newcommand*{\RR}{\mathbb{R}}
\newcommand*{\ZZ}{\mathbb{Z}}

\renewcommand*{\vec}{\mathbf}

\begin{document}

\title{Pulling Subdivisions of Cones and Blowups of Monomial Ideals on Affine Toric Varieties}
\date{\today}
\author{Howard M Thompson}
\address{Mathematics Department\\ 402 Murchie Science Building \\ 303 East Kearsley Street \\ Flint, MI 48502-1950}
\email{\href{mailto:hmthomps@umflint.edu}{hmthomps@umflint.edu}}
\urladdr{\url{http://homepages.umflint.edu/~hmthomps/}}
\thanks{}
\subjclass[2010]{Primary: 14M25; Secondary: 52B20}
\keywords{Multiplier ideals, toric varieties}
\begin{abstract}
  This short note solves the following problem: Given a map $\pi:X_\Sigma\rightarrow X_\sigma$ of normal toric varieties corresponding to a coherent subdivision $\Sigma$ of a cone $\sigma$, find an ideal sheaf $\sheaf{I}$ on $X_\sigma$ such that $\pi$ is the blowup of $\sheaf{I}$.
\end{abstract}
\maketitle

\section{Introduction}

Thompson~\cite{MR3448804} takes advantage of the following theorem to find identify a well-controlled log resolution of a monomial space curve.

\begin{thm}(Gonz{\'a}lez P{\'e}rez \& Teissier~\cite[Theorem~3.1]{MR1892938}) \label{t:GT}
  Let $\sigma\subset N_\RR$ be strictly convex rational polyhedral cone, let $\phi:\monoid_\sigma\rightarrow\monoid$ be a surjective homomorphism of pointed affine semigroups, let $\phi_\RR:M_\RR\rightarrow\RR\monoid$ be the induced linear map, let $\ell=\ker(\phi_\RR)^\perp\subset N_\RR$, let $\Bbbk$ be a field, and let $Z=\spec\left(\Bbbk[\monoid]\right)\subset X_\sigma$. If $\Sigma$ is any subdivision of $\sigma$ containing the cone $\tau=\sigma\cap\ell$ and $N_\tau=\ZZ(N\cap\tau)$, then
  \begin{enumerate}[(1)]
    \item  The strict transform of $Z$ by the morphism induced by the subdivision $\pi_\Sigma:X_\Sigma\rightarrow X_\sigma$ is contained in $X_\tau$, it is isomorphic to $X_{\tau,N_\tau}$ and the restriction $\pi_\Sigma|_{X_{\tau,N_\tau}}:X_{\tau,N_\tau}\rightarrow Z$  is the normalization map.
    \item The morphism $\pi_\Sigma$ is a partial embedded resolution of $Z\subset X_\sigma$. (That is, any toric desingularization of $X_\Sigma$ provides an embedded resolution of $Z\subset X_\sigma$).
  \end{enumerate}
\end{thm}

In \vref{s:pull}, we  will construct a minimal subdivision as in the theorem by defining pulling subdivision for cones and present a simple example. In \vref{s:ideals}, we briefly describe how to produce an ideal from a coherent subdivision.

Our general reference on toric varieties is Cox, Little \& Schenck~\cite{MR2810322} and our general reference on subdivisions is De~Loera, Rambau \& Santos~\cite{MR2743368}.

\section{Pulling Subdivision} \label{s:pull}

Let $\sigma,\tau\subset N_\RR$ be strictly convex rational polyhedral cones such that $\tau$ is a subset of $\sigma$. We do not assume $\tau$ is a subcone of $\sigma$.

\begin{defn}
  The following construction defines the \emph{pulling subdivision} $\pull_\tau\sigma$ of $\tau$ in $\sigma$. Fix a (rational) hyperplane $H$ not containing the origin such that $H\cap\rho$ is nonempty for each ray $\rho\in\sigma(1)$. Let
  \[
    \mathscr{A}=\{H\cap\rho\mid\rho\in\sigma(1)\}\text{ and }\mathscr{B}=\{H\cap\rho\mid\rho\in\tau(1)\}.
  \]
  be (necessarily finite) point sets in $N_\RR\times\RR$, Consider the polytope $\convhull((\mathscr{A}\times\{0\})\cup(\mathscr{B}\times\{1\}))$ corresponding to the height function $\omega:\mathscr{A}\rightarrow\RR$ that is $1$ on $\mathscr{B}$ and $0$ on $\mathscr{A}\setminus\mathscr{B}$. By construction, the projection of the upper hull of this polytope onto the first factor is a coherent subdivision of $P=H\cap\sigma$. This subdivision of $H\cap\sigma$, $\pull_Q P$, is (essentially) the result of pulling $Q=H\cap\tau$ as in Section~2.2 of Haase \& Zharkov~\cite{HZ2002}. Let $\pull_\tau\sigma$ be the fan consisting of the cones over the faces of $\pull_Q P$. 
\end{defn}

Note that the height function $\omega:\mathscr{A}\rightarrow\RR$ that is $1$ on $\mathscr{B}$ and $0$ on $\mathscr{A}\setminus\mathscr{B}$ extends to a support function $\varphi:|\pull_\tau\sigma|\rightarrow\RR$ given by rational vectors $\{\vec{u}_\rho\}_{\rho\in\pull_\tau\sigma}\subset\QQ M$. So, some multiple of $\varphi$ is a support function that is integral with respect to the lattice $N$. Let this multiple be given by the set $\{\vec{m}_\rho\}_{\rho\in\pull_\tau\sigma}$. This set is the Cartier data for some Cartier divisor as in Cox, Little \& Schenck~\cite[Theorem~4.2.8]{MR2810322}.

\begin{ex}
  Let $\vec{n}_1=\begin{bmatrix} 1 & 0 & 0 \end{bmatrix}^\mathsf{T}$, let $\vec{n}_2=\begin{bmatrix} 0 & 1 & 0 \end{bmatrix}^\mathsf{T}$, let $\vec{n}_3=\begin{bmatrix}0 & 0 & 1 \end{bmatrix}^\mathsf{T}$, let $\vec{n}_4=\begin{bmatrix} 2 & 1 & 0 \end{bmatrix}^\mathsf{T}$, let $\vec{n}_5=\begin{bmatrix}0 & 1 & 2 \end{bmatrix}^\mathsf{T}$, let $\sigma=\RR_{\geq0}\vec{n}_1+\RR_{\geq0}\vec{n}_2+\RR_{\geq0}\vec{n}_3$ be the positive orthant in $\RR^3$, let $\tau=\RR_{\geq0}\vec{n}_4+\RR_{\geq0}\vec{n}_5$, and pick $\begin{bmatrix} 1 & 1 & 1 \end{bmatrix}v=1$ for $H$. Then, $\mathscr{A}=\{\vec{n}_1,\vec{n}_2,\vec{n}_3\}$ and $\mathscr{B}=\left\{\frac{1}{3}\vec{n}_4,\frac{1}{3}\vec{n}_5\right\}$.

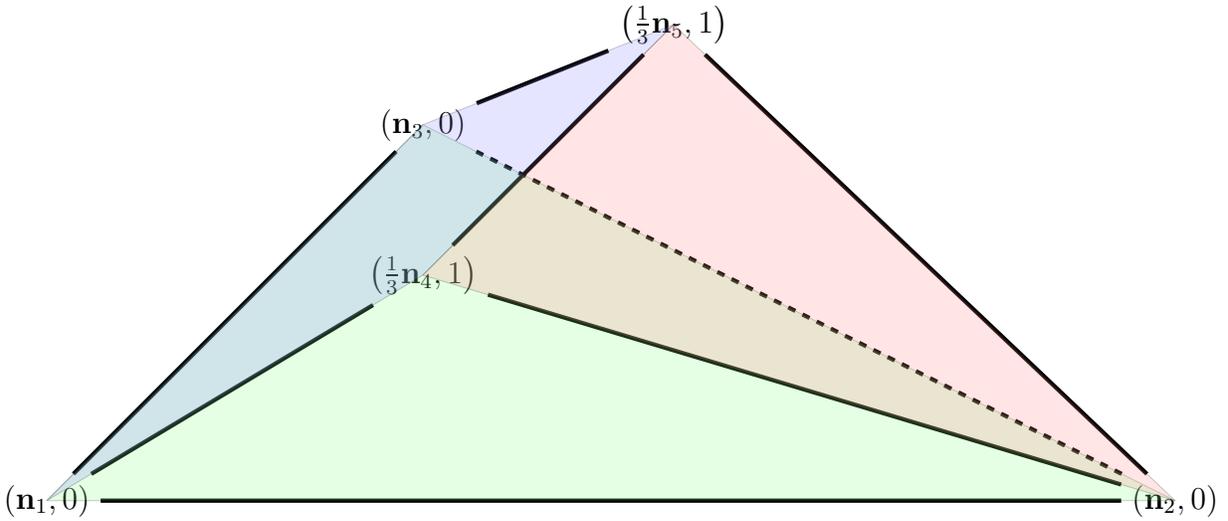
  \begin{center}
  \begin{tikzpicture}
    \node (N1) at (0,0) {$(\vec{n}_1,0)$};
    \node (N2) at (15,0) {$(\vec{n}_2,0)$};
    \node (N3) at (5,5) {$(\vec{n}_3,0)$};
    \node (V4) at (5,3) {$\left(\frac{1}{3}\vec{n}_4,1\right)$};
    \node (V5) at (8.333,6.333) {$\left(\frac{1}{3}\vec{n}_5,1\right)$};
    \draw [ultra thick] (N1) -- (N2);
    \draw [ultra thick] (N1) -- (N3);
    \draw [ultra thick,dashed] (N2) -- (N3);
    \draw [ultra thick] (N1) -- (V4);
    \draw [ultra thick] (V4) -- (N2);
    \draw [ultra thick] (N2) -- (V5);
    \draw [ultra thick] (V5) -- (N3);
    \draw [ultra thick] (V4) -- (V5);
    \draw [fill=green!50,opacity=0.2] (0,0) -- (15,0) -- (5,5) -- (0,0);
    \draw [fill=blue!50,opacity=0.2] (0,0) -- (5,5) -- (8.333,6.333) -- (5,3) -- (0,0);
    \draw [fill=red!50,opacity=0.2] (15,0) -- (5,3) -- (8.333,6.333) -- (15,0);
  \end{tikzpicture}
  \end{center}
  
  In the figure, the lower facet $\convhull(\mathscr{A}\times\{0\})$ is green, the upper facets containing $\convhull(\mathscr{B}\times\{1\})$ are blue and red, and the other two (unshaded) facets are vertical. The projections of the two vertical facets of the polyhedron to the green facet are not maximal cells of the subdivision. The projections of the blue and red facets onto the green facets yield the maximal cells of $\pull_Q P$.

  According to a Macaulay2~\cite{M2} calculation, the $3$-dimensional polytope in $\RR^4$ with vertex set $(\mathscr{A}\times\{0\})\cup(\mathscr{B}\times\{1\})$ is given by the halfspaces
  \[
    \begin{matrix*}[r]
        w_1 & + w_2 & -3w_3 &       & \leq 1 &     \\
      -3w_1 & + w_2 &  +w_3 &       & \leq 1 &     \\
            &       &       &  -w_4 & \leq 0 &     \\
      -3w_1 & +3w_2 & -3w_3 & +4w_4 & \leq 3 & \qquad (1) \\
       3w_1 & -9w_2 & +3w_3 & +4w_4 & \leq 3 & \qquad (2)
    \end{matrix*}
  \]

  The triangle with vertex set $\left\{\vec{n}_2,\frac{1}{3}\vec{n}_4,\frac{1}{3}\vec{n}_5\right\}$ is the projection of the red facet supported on the hyperplane defining the halfspace $(1)$ and the trapezoid with vertex set $\left\{\vec{n}_1,\vec{n}_3,\frac{1}{3}\vec{n}_4,\frac{1}{3}\vec{n}_5\right\}$ is the projection of the blue facet supported on the hyperplane defining the halfspace $(2)$. Thus, the maximal cones of the fan $\pull_\tau\sigma$ are $\sigma_1=\RR_{\geq0}\vec{n}_2+\RR_{\geq0}\vec{n}_4+\RR_{\geq0}\vec{n}_5$ and $\sigma_2=\RR_{\geq0}\vec{n}_1+\RR_{\geq0}\vec{n}_3+\RR_{\geq0}\vec{n}_4+\RR_{\geq0}\vec{n}_5$. Note that this fan is not refined by nor a refinement of the star subdivision along either of the rays of  $\tau$. Here the support function $\varphi$ extending the height function defining $\pull_\tau\sigma$ is given by $\vec{u}_{\sigma_1}=\begin{bmatrix} \frac{3}{2} & 0 & \frac{3}{2} \end{bmatrix}$ and $\vec{u}_{\sigma_2}=\begin{bmatrix} 0 & 3 & 0 \end{bmatrix}$. That is,
  \[
    \varphi:|\pull_\tau\sigma|\rightarrow\RR; \qquad \vec{v}\mapsto\min\{\langle\vec{u}_1,\vec{v}\rangle,\langle\vec{u}_2,\vec{v}\rangle\}.
  \]
  So, $2\varphi$ is given by $\vec{m}_{\sigma_1}=\begin{bmatrix} 3 & 0 & 3 \end{bmatrix}$ and $\vec{m}_{\sigma_2}=\begin{bmatrix} 0 & 6 & 0 \end{bmatrix}$. Here, $\pull_\tau\sigma$ is the (inward) normal fan of the polyhedron $\convhull((\vec{m}_{\sigma_1}+\monoid_\sigma)\cup(\vec{m}_{\sigma_2}+\monoid_\sigma))$. Let $x_i=\chi^{\vec{e}_i}$ where $\vec{e}_i$ is the $i$th standard basis element of $M_\RR$ and let $I=\overline{\left(x_1^3x_3^3,x_2^6\right)}$ in $\Bbbk[\monoid_\sigma]=\Bbbk[x_1,x_2,x_3]$. In this case, $\pull_\tau\sigma$ is the fan of the blowup of $I$. 
\end{ex}

  More generally, if $\newt(I)=\convhull\{\vec{m}\mid\chi^{\vec{m}}\in I\}$ is the Newton polyhedron of a torus invariant ideal $I\subset\Bbbk[\monoid_\sigma]$ for some strictly convex rational polyhedral cone $\sigma$, then the dual fan of $\newt(I)$ is the fan of the normalized blowup of $I$. See Thompson~\cite{hT2003}.

\section{From Coherent Subdivisions to Ideals} \label{s:ideals}

Now, if $\Sigma$ is a fan that is a coherent subdivision of a fan $\Delta$ in $N_\RR$, then we would like to find an ideal sheaf $\sheaf{I}$ on $X_\Delta$ such that the map of toric varieties $\pi_\Sigma:X_\Sigma\rightarrow X_\Delta$ induced by the subdivision is the blowup of $\sheaf{I}$. To achieve this, we want a support function $\varphi$ that is integral with respect to the lattice $N$ with Cartier data $\{\vec{m}_\sigma\}_{\sigma\in\Sigma}$ (as in Cox, Little \& Schenck~\cite[Theorem~4.2.8]{MR2810322}) such that whenever $\tau\in\Delta$ with $\sigma\subset\tau$, $\vec{m}_\sigma\in\monoid_\tau$. Solve the system of linear inequalities for the subset $\{\vec{m}_\sigma\}_{\sigma\in\Sigma_{max}}$. Namely,
\begin{enumerate}
  \item $\langle\vec{m}_{\sigma_1},\vec{v}_\rho\rangle=\langle\vec{m}_{\sigma_2},\vec{v}_\rho\rangle$ whenever $\vec{v}_\rho$ is the primitive vector on a ray $\rho\subset(\sigma_1\cap\sigma_2)$;
  \item $\langle\vec{m}_{\sigma_1},\vec{v}_\rho\rangle<\langle\vec{m}_{\sigma_2},\vec{v}_\rho\rangle$ whenever $\vec{v}_\rho$ is the primitive vector on a ray $\rho\subset(\sigma_1\setminus\sigma_2)$.
\end{enumerate}
Choose the torus invariant sheaf $\sheaf{I}$ on $X_\Delta$ such that, for each $\tau\in\Delta_{max}$,
\[
  \Gamma(U_\tau,\sheaf{I})=\overline{(\chi^{\vec{m}_\sigma}\mid\sigma\subset\tau)}\subset\Bbbk[\monoid_\tau]
\]

\end{document}